\documentclass[12pt]{article} 
\usepackage[utf8x]{inputenc}
\usepackage{tikz} 
\usepackage[all]{xy}
\usepackage{authblk}
\usepackage[babel=true]{csquotes}
\usepackage{amsmath,amsthm, amssymb,amsfonts}
\usepackage[hmargin=1in,vmargin=1in]{geometry}
\numberwithin{equation}{subsection} 
\xyoption{arc}
\usepackage{eucal}
\usepackage{indentfirst}
\usepackage{mathrsfs}
\usepackage{verbatim}
\usepackage{makeidx}
\usepackage{graphicx}
\usepackage{yfonts} 
\usepackage{hyperref}
\usepackage{enumitem} 
\usepackage{extpfeil} 
\usepackage{dsfont}
\usepackage[T1]{fontenc}
\newtheorem{thm}{Theorem}[section]

\newtheorem{conj}[thm]{Conjecture}

\newtheoremstyle{bidule}
{3pt}
{3pt}
{}
{}
{\scshape}
{.}
{.5em}
{}

\theoremstyle{definition}
\newtheorem{rmk}[thm]{Remark}

\newcommand{\C}{\mathcal{C}}
\newcommand{\Ub}{\mathcal{U}}

\newcommand{\F}{\mathcal{F}}

\newcommand{\Z}{\mathbb{Z}}
\newcommand{\R}{\mathbb{R}}

\newcommand{\Cx}{\mathbb{C}}

\renewcommand{\to}{\longrightarrow}

\newcommand{\N}{\mathbb{N}}

\renewcommand{\k}{\textbf{k}} 

\newcommand{\tx}{\text}

\renewcommand{\to}{\longrightarrow}
\DeclareMathOperator\Id{Id}
\DeclareMathOperator\Hom{Hom}

\DeclareMathOperator\Set{\textbf{Set}} 
\DeclareMathOperator\Cat{\mathbf{Cat}}

\DeclareMathOperator\op{\tx{op}}

\DeclareMathOperator\Ho{\mathbf{ho}} 

\DeclareMathOperator\Vect{\textbf{Vect}} 

\DeclareMathOperator\sset{\mathbf{sSet}} 
\DeclareMathOperator\ssetj{\mathbf{sSet}_{J}} 
\DeclareMathOperator\ssetq{\mathbf{sSet}_{Q}} 
 
\DeclareMathOperator\Top{\mathbf{Top}}
\DeclareMathOperator\Var{\mathbf{Var}}
\DeclareMathOperator\Mod{\mathbf{Mod}}

\title{Symmetries} 
\author{Hugo V. Bacard \thanks{\textit{E-mail address}: \href{mailto:hbacard@uwo.ca}{hbacard@uwo.ca}
}}
 \affil{Western University}
\date{July 5, 2014}
\begin{document}
\maketitle
\begin{abstract}
Classical mathematics are founded within set theory, but sets don't have \emph{symmetries}. We conjecture that if we allow sets with symmetries, then many problems such as \emph{Mirror symmetry} or \emph{Homological mirror symmetry} can be explained. One way to do this is to embed sets to higher categories and especially into higher groupoids as already envisioned by Grothendieck,  and more recently by Voevodsky. We simply outline this idea in these notes.
\end{abstract}
\setcounter{tocdepth}{1}
\tableofcontents

\vspace*{2cm}

\hfill « \textit{Tout est intéressant !}»

\hfill Carlos Simpson

\vspace*{2cm}

\section{Symmetries in maths}
Given a set $X= \{ a,b,c,..\}$ such as the natural numbers $\N= \{0,1,...,p,... \}$, there is a standard procedure that amounts to regard $X$ as a category with only identity morphisms. This is the \emph{discrete functor} that takes $X$ to the category denoted by $Disc(X)$ where the hom-sets are given by $\Hom(a,b)= \emptyset$ if $a\neq b$, and $\Hom(a,b)= \{Id_a\}= 1$ if $a=b$. $Disc(X)$ is in fact a groupoid.\\

But in category theory, there is also a procedure called \emph{opposite or dual}, that takes a general category $\C$ to its opposite $\C^{op}$. And to put the reader on the road, let's  also call $\C^{op}$ the \emph{reflection of $\C$} by the \emph{mirror functor} $(-)^{op}$.   

Now the problem is that if we restrict this procedure to categories such as $Disc(X)$, there is no way to distinguish $Disc(X)$ from $Disc(X)^{op}$. And this is what we mean by \emph{sets don't show symmetries}. In the program of Voevodsky, we can interpret this by saying that: 
\begin{center} 
\emph{`The identity type is not good for sets, instead we should use the Equivalence type. But to get this, we need to move to from sets to Kan complexes i.e., $\infty$-groupoids'}.
\end{center}

So far we've used set theory with this lack of symmetries, as foundations for mathematics. And some \emph{symmetry phenomenons} occur as we progress in maths, and sometimes we're unable to figure out why exactly.\\

Grothendieck \cite{SGA1} has already seen this when he moved from sheaves of sets, to sheaves of groupoid (stacks), because he wanted to allow objects to have symmetries (automorphisms). If we look at the Giraud-Grothendieck  picture on nonabelian cohomology \cite{Gir}, then what happens is an extension of coefficients $\Ub:\Set \hookrightarrow \Cat$.  But this embedding is too big as we mentioned in \cite{Bacard_Und_QS}. Rather we should consider first the comma category $\Cat\downarrow \Ub$, whose objects are functors $C \to Disc(X)$. And then we should consider the full subcategory consisting of  functors $C \xrightarrow{\sim} Disc(X)$ that are equivalences of categories. This will force $C$ to be a groupoid, that looks like a set. And we call such $C \xrightarrow{\sim} Disc(X)$ a \emph{Quillen-Segal } $\Ub$-object.\\

This category of Quillen-Segal objects should be called the \emph{category of sets with symmetries}. Following Grothendieck's point of view, we've denoted by $\Cat_\Ub[\Set]$ the comma category, and think of it as \emph{categories with coefficients or coordinates in sets}. This terminology is justified by  the fact that the functor $\Ub:\Set \hookrightarrow \Cat$ is a morphism of (higher) topos, that defines a \emph{geometric point in $\Cat$}. The category of set with symmetries is like the homotopy neighborhood of this point, similar to a one-point going to a disc or any contractible object. The advantage of the Quillen-Segal formalism is  the presence of a Quillen model structure on $\Cat_\Ub[\Set]$ such that the fibrant objects are Quillen-Segal objects (\cite[Theorem 8.2]{Bacard_QS}, \cite[Theorem 1.2 ]{Bacard_Und_QS}. 

In standard terminology this means that if we embed a set $X$ in $\Cat$ as $Disc(X)$,  and take an `projective resolution' of it, then we get an equivalence of groupoids $P \xrightarrow{\sim} Disc(X)$, and $P$ has symmetries. Concretely what happens is just a factorization of the identity (type) $\Id: Disc(X) \to Disc(X)$ as a cofibration followed by a trivial fibration:

$$Disc(X) \hookrightarrow P \xrightarrow{\sim} Disc(X).$$ 
The first morphism is also automatically an equivalence.\\

We regard in \cite{Bacard_Und_QS} this process of embedding $\Set \hookrightarrow QS\{\Cat_\Ub[\Set] \}$ as \emph{a minimal homotopy enhancement}. The idea is that there is no good notion of homotopy (weak equivalence)  in $\Set$, but there are at least two notions in $\Cat$: equivalences of categories and the   equivalences of classifying spaces \emph{à la} Grothendieck-Kan-Quillen-Segal-Thomason.\\

This last class of weak equivalences is important for what we believe happens with \emph{mirror phenomenons}. We isolated the discussion in the next paragraph. But for experts: the mirror of a manifold should be the \emph{opposite of its fundamental Poincaré $\infty$-groupoid}. We make a precise statement below.

\section{Mirrors} 

Given a compact K\"ahler manifold $Y$, we know that the cohomology groups $H^\star(Y,\Cx)$ have a Hodge decomposition $H^{p,q}$ (see \cite{Hodge_3} \cite{Voisin_livre_1,Voisin_livre_2}). Now because we have Poincaré duality, and the comparisons between singular and De Rham cohomologie, we know that any other space $Z$ that has the same homotopy type as $Y$ will have the same cohomology groups. Consequently they will share the same \emph{Hodge diamond}, thus its symmetries.\\

This means that the symmetry of the Hodge diamond is mostly attached to the homotopy type of $Y$. This is not surprising anymore because it can already be seen from the equivalence of the  De Rham cohomology which is analytic, and the Betti cohomology which is something purely simplicial. In fact, it can also be seen from the (smooth) homotopy invariance of  De  Rham cohomology.\\

We believe that this symmetry can be understood using the Quillen-Segal formalism as follows. Given $Y$, let's consider $Y_{top} \in \Top$. Recall that we have a Quillen equivalence $\Ub:Top \to \ssetq$, where $\Ub=Sing$ is the singular functor whose left adjoint is the geometric realization. When we consider the comma category $\ssetq_\Ub[Top]= \sset\downarrow \Ub$, we are literally creating in French a \emph{``trait d'union''}, between the two categories. And when we consider the subcategory of Quillen-Segal objects, then our result \cite[Theorem 1.2 ]{Bacard_Und_QS} says that we have a triangle that descends to a triangle of equivalences between the  homotopy categories. In fact there is a much better statement.

\[
\xy
(0,18)*+{(\ssetq \downarrow \Ub)}="W";
(0,0)*+{\Top}="X";
(30,0)*+{\ssetq}="Y";
{\ar@{->}^-{\Ub}"X";"Y"};
{\ar@{->}^-{}"W";"X"};
{\ar@{->}^-{s}"W";"Y"};
\endxy
\] 

It turns out that if we apply \cite[Theorem 1.2 ]{Bacard_Und_QS} to the same functor but we choose the Joyal model structure $\ssetj$, we get the \emph{Homotopy hypothesis} (see\cite{Ara_gpds} for the statement of this hypothesis).\\

A  fibrant replacement  of $Y$ in the model category $\ssetq_\Ub[\Top]$, is a trivial fibration $\F \xtwoheadrightarrow{\sim} \Ub(Y)$, where $\F$ is fibrant in $\ssetq$, that is a Kan complex. But a Kan complex is exactly an $\infty$-groupoid. $\infty$-Groupoids generalize groupoids, and still are category-like. In particular we can take their opposite (or dual), just like we consider the opposite category $\C^{op}$ of a usual category, as outlined in the beginning. 

\begin{conj}
Given $Y$ as above, we can think of the mirror of $Y$ as the opposite $\infty$-groupoid $\F^{op}$. A good approximation of $\F^{op}$ can be obtained by the \emph{schematization functor} à la To\"en applied to  the simplicial set (quasicategory) underlying $\F^{op}$. \\

We can take as model for $\F$ the fundamental $\infty$-groupoid $\Pi_{\infty}(Y)$. And depending on the dimension it's enough to stop at the corresponding $n$-groupoid.\\

To\"en schematization functor can also be obtained from the Quillen-Segal formalism applied to the embedding 
$$\Ub: Sh(\Var(\Cx)) \hookrightarrow  sPresh(\Var(\Cx),$$
where on the right hand side we consider the model category of simplicial presheaves à la Jardine-Joyal. The representability of  the $\pi_0$ of the schematization has to be determined by descent along the equivalence type. 
\end{conj}

\section{Conjectures}

We now list some conjectures. We use the same notations as in \cite{Bacard_Und_QS}. These statements are only inspired by some abstract thinking and not by experience in algebraic geometry.

\begin{conj}
\begin{enumerate}
\item As mentioned above, we can enhance $Top \hookrightarrow \ssetj[Top]$ by looking at the Quillen-Segal objects. Then given $Y$ as before, then the mirror of $Y$ should correspond to thhe\footnote{'thhe' is the homotopy version of 'the' (Drinfeld)} opposite $\infty$-groupoid $\Pi_\infty(Y)$. Indeed by a theorem of To\"en \cite{Toen_axiom_hcat}, we know that $\Z/2= Gal(\Cx/\R)$ acts on $\Ho(\ssetj)$, where $0$ is the identity and $1$ is the opposite-category construction. In particular we have $\pi_1(Y^{op})= \pi_1(Y)^{op}$, as expected. \\
 It seems that it's not surprising that $\Z/2$ appears in supersymmetry.

\item The generating trivial cofibration in the folk model structure on $\Cat$ is the minimal equivalence of groupoids 
$$\ast \hookrightarrow \{ \underbrace{\ast \leftrightarrows \ast}_{=\tx{walking iso}} \}.$$
We believe that this should be considered as the introduction of the  \emph{Higgs boson}, and therefore we shall call it the \emph{Higgs} equivalence or \emph{Higgs symmetry}. We are tempted to denote the two arrows in the interval category $\{\ast \leftrightarrows \ast \}$ by $e^+$ and $e^-$. 

\item We have a factorization of the identity $\Id_\ast$ through the Higgs boson, as a Feynman-like diagram:

\begin{equation}\label{boson}
\ast \xrightarrow{\Id} \ast= \ast \hookrightarrow \{\ast_ \leftrightarrows \ast \} \to \ast.
\end{equation}

We would like to interpret this as the boson hiding between the identity, and therefore invisible. 
\item We believe that the fact that $\Cx$ is algebraically closed and that $Gal(\Cx/\R)= \Z/2$ is \emph{quantumly related} to Higgs boson. The reason being  that following the factorization \eqref{boson}, we are tempted to write:
$$e^+ e^- \simeq \Id \quad \leadsto \quad i (-i)= 1.$$ 
It seems that this phenomenon explains why representation of fundamental group of complex projective variety leads to Higgs bundle as shown in the work of Hitchin \cite{Hitchin_higgs} and Simpson \cite{Simpson_Rep_II,Simpson_Rep_I}.

\item We can enhance $\Vect_\Cx \hookrightarrow (\infty,n)\Cat[\Vect_\Cx]$ with the relative pushout product. Taking $(\infty,n)\Cat[\Vect_\Cx]$ as coefficient for TQFT, should explain why TQFT are classified by fully dualizable objects. Indeed To\"en's theorem  has been generalized by Barwick and Schommer-Pries \cite{Barwick_SPries}. There is an action of  $(\Z/2)^n$ on $\Ho[(\infty,n)\Cat]$, that correspond to the different opposite (=mirror) constructions for $1$-morphisms, to $n$-morphisms. This action should explain the \emph{$n$-dualizable objects}, as named by Lurie.

\item Similarly we can enhance $n\tx{-Fold} \hookrightarrow (\infty,n)\Cat[n\tx{-Fold}]$ to see the symmetries at every level for manifolds. Taking $n=4$ for space-time should fix some issues occurring in Physics with set theory.

\item We can enhance $dg\tx{-}\Cat \hookrightarrow (\infty,2)\Cat[dg\tx{-}\Cat]$ and similarly we have an action of $(\Z/2)^2$ on the homotopy category $\Ho[(\infty,2)\Cat]$. This action should explain the proof of Deligne's conjecture given by Tamarkin \cite{Tamarkin_Deligne}, in particular why there are $2$-discs acting on the Hochschild cohomology. The number $n=2$ is the exponent in $(\Z/2)^n$. It should also agree with Tamarkin's answer to Drinfeld's question:
\begin{center}
\emph{What do dg-categories form ?}
\end{center} 

This philosophy should fit in Kontsevich's program on homological mirror symmetry.
\item We can enhance $\Mod_\k \hookrightarrow (\infty,1)\Cat[\Mod_\k]$ with the relative pushout product. It should be interesting to let cohomology theories in algebraic geometry take their coefficient in this enhancement. 
\end{enumerate}
\end{conj}

\begin{rmk}

\begin{enumerate}
\item If we follow our philosophy, it's not surprising that there is no direct link between a variety $Y$ and its mirror $Y^{\op}$. Because there is no direct link in general between a classical category $\C$ and its opposite $\C^{op}$, unless $\C$ is Tannakian or at least has duals.

\item This last fact should have its analogy with duality in geometry such as Poincaré, Serre, Grothendieck dualities. Lurie has already taken this direction when he speaks of \emph{Poincaré object} and \emph{nonabelian Poincaré duality}. Simpson has further developed his program on nonabelian Hodge theory using higher stacks.  And there is a long list of people who are currently developing these ideas that aim to fix the issues caused by set theory and its lack of symmetries. 

\item It would be interesting to understand the statement of the Hodge conjecture in terms of $\infty$-groupoids that are fixed by the homotopy action of $\Z/2= Gal(\Cx/\R)$ and its power $(\Z/2)^n$ . After all, as mentioned before, the cohomology of a subvariety is actually a homotopy theory invariant.
\end{enumerate}
\end{rmk}

\vspace*{1cm}
\bibliographystyle{plain}
\bibliography{Bibliography_LP_COSEG}
\end{document}